\documentclass[11pt]{article}
\usepackage[T1]{fontenc}
\usepackage[toc,page]{appendix} 
\usepackage{amsmath}
\usepackage{amsfonts}
\usepackage{amssymb}
\usepackage{geometry}
\usepackage{fancyhdr}
\author{Vincent ~NOLOT\\ \vspace{5pt} \\
   Institut de Mathématiques de Bourgogne,\\
   Université de Bourgogne, 21078 Dijon, France.\\
   \texttt{vincent.nolot@u-bourgogne.fr}}
\title{Optimal transport on the classical Wiener space with different norms}

\begin{document}

\maketitle
\newtheorem{thm}{\textsc{Theorem}}[section]
\newtheorem{lemme}{\textsc{Lemma}}[section]
\newtheorem{prop}{\textsc{Proposition}}[section]
\newtheorem{cor}{\textsc{Corollary}}[section]
\newtheorem{defin}{\textsf{Definition}}[section]
\newtheorem{ex}{\textsf{Exemple}}[section]
\newenvironment{rmq}{\noindent \underline{Remark:}}{}
\newenvironment{nota}{\textsf{Notations:}}{}
\newenvironment{changemargin}[2]{\begin{list}{}{%
\setlength{\topsep}{0pt}%
\setlength{\leftmargin}{0pt}%
\setlength{\rightmargin}{0pt}%
\setlength{\listparindent}{\parindent}%
\setlength{\itemindent}{\parindent}%
\setlength{\parsep}{0pt plus 1pt}%
\addtolength{\leftmargin}{#1}%
\addtolength{\rightmargin}{#2}%
}\item }{\end{list}}
\newenvironment{demo}{\noindent \texttt{\textit{Proof:}\\}}{\rule{1ex}{1ex}\vspace{15pt}}
\newenvironment{details}
{\noindent $\rightarrow$ \hspace{15pt}}
{}
\renewcommand{\theequation}{\thesection.\arabic{equation}}

\pagestyle{fancy}
\renewcommand{\sectionmark}[1]{\markboth{#1}{}}
\renewcommand{\headrulewidth}{0.5pt}
\renewcommand{\footrulewidth}{0.5pt}

\abstract{In this paper we study two basic facts of optimal transportation on Wiener space $W$. Our first
aim is to answer to the Monge Problem on the Wiener space endowed with the Sobolev type norm
$\|.\|_{k,\gamma}^p$ with $p\geq 1$ (cases $p=1$ and $p>1$ are considered apart). 
The second one is to prove $1-$convexity (resp. $1/C_{k,\gamma}^2-$convexity) along (constant speed)
geodesics of relative entropy in $(\mathcal{P}_2(W),W_2)$, where $W$ is endowed
with the infinite norm (resp. with $\|.\|_{k,\gamma}$), and $W_2$ is the $2-$distance of Wasserstein.}


\section{Introduction}
We are interested in two problems in optimal transportation on Wiener space.
We refer to a recent work of Ambrosio and Gigli
\cite{USERG} for a survey and basic tools of optimal transport theory.\\

At first we answer to the Monge Problem on the Wiener space relatively to the cost (introduced by Airault and Malliavin
in \cite{AIRMALL})
$\|.\|_{k,\gamma}^p$ defined as
$$\|w\|_{k,\gamma}:=\left(\int_0^1\int_0^1\frac{(w(t)-w(s))^{2k}}{|t-s|^{1+2k\gamma}}dtds\right)^{1/2k},$$
for suitable parameters $k,\gamma$.

Monge Problem is largely considered in several settings since for years. 
Nowadays there are a lot of means to prove the existence of an optimal map resolving this Problem,
which are summarized in \cite{AMB0}. 

Recently, there was considerable advances concerning Monge Problem in $\mathbb{R}^n$.
First in 1996, Gangbo and McCann solved Monge Problem in \cite{GANGMC} when the cost is strictly convex.
Then people are interesting in the case of different norms on $\mathbb{R}^n$. Indeed the problem
becomes more difficult since a norm is never strictly convex.
In 2003, Monge Problem was solved when the cost was a crystalline norm
by Ambrosio, Kirchheim and Pratelli in \cite{AMB2}. When the cost is a general norm,
Monge Problem was solved independently by Champion, De Pascale in \cite{CHAMP} and by
Caravenna in \cite{CARAV} in 2010. All of these latter cases we lose unicity of optimal map.

In this paper, we turn our attention on the Wiener space, an infinite dimensional space.
The Monge Problem in the Wiener space was solved by Feyel and Ustunel in \cite{FEY}
for the cost $|.|_H^2$ induced by Cameron-Martin norm $|.|_H$ which is Hilbertian. With this cost,
we have unicity of optimal map. Strategy of latter authors was to pass by finite dimensional 
approximations of Wiener space: the aim being to reduce
the Monge Problem on Wiener space onto finite dimensional spaces and applying known results. Then they used a
selection theorem to go back up on Wiener space. 
Then Monge Problem was solved by Cavalletti in \cite{CAVA} for the cost $|.|_H$, again passing by
finite dimensional approximations.

We are interested to endow the Wiener space $(W,H,\mu)$ (where $\mu$ is the Wiener measure) with two other natural norms:
the infinite norm $|.|_\infty$ and the Sobolev type norm $\|.\|_{k,\gamma}$. For the first one,
Monge Problem is still open, and we can expect not to have unicity of optimal map, providing
it exists somehow. Here the norm considered $\|.\|_{k,\gamma}$ is not Hilbertian
and is weaker than the Cameron-Martin norm in sense that for some $C_{k,\gamma}>0$:
\begin{equation}\label{inqNorms}
	\|x\|_{k,\gamma}\leq C_{k,\gamma} |x|_H~~~~\mathrm{for~all}~x\in W.
\end{equation}
Let us emphasize that the right hand side is equal to infinite $\mu-$almost everywhere, because of zero measure of
Cameron-Martin space. Nevertheless our norm $\|.\|_{k,\gamma}$ satisfies suitable conditions presented in section 2.1.
The first aim of our paper is to solve Monge Problem for the cost
$\|.\|_{k,\gamma}^p$:
\begin{eqnarray}\label{MongeP}
	\inf_{G_\# \rho_0=\rho_1} \int_{W\times W}\|G(x)-x\|_{k,\gamma}^pd\rho_0(x).
\end{eqnarray}

In other words we will establish following theorem:

\begin{thm}
Let $\rho_0$ and $\rho_1$ be two measures on $W$ satisfying condition (\ref{PreservingMass}) below.
\begin{enumerate}
	\item If $p>1$ and $\rho_0$ is absolutely continuous with respect to $\mu$, then
	there exists a map $T$ unique up to a set of zero measure for $\rho_0$, which minimizes (\ref{MongeP}).
	Moreover there is a unique optimal transference plan between $\rho_0$ and $\rho_1$ relatively to the cost $\|.\|_{k,\gamma}^p$,
	which is exactly $(Id\times T)_\#\rho_0$.
	\item If $p=1$ and both $\rho_0$ and $\rho_1$ are absolutely continuous with respect to $\mu$, then
	there exists an optimal transference plan $\Pi$ between $\rho_0$ and $\rho_1$ relatively to the cost $\|.\|_{k,\gamma}$,
	such that $\Pi$ is concentrated on a graph of some map $T$ minimizing (\ref{MongeP}).
\end{enumerate}
\end{thm}

The second purpose of this paper is to find a lower bound of Ricci curvature in the sense
of Sturm in \cite{STURM}, for Wiener space endowed with the infinite norm $|.|_\infty$ or the norm $\|.\|_{k,\gamma}$.
 It is relied
to weak $K-$convexity along geodesic of relative entropy. More precisely we will prove:

\begin{thm}
If $\rho_0$ and $\rho_1$ are probability measures on $W$ both absolutely continuous with respect
to $\mu$, then there exists some (constant speed) geodesic $\rho_t$ induced by an optimal transference between $\rho_0$
and $\rho_1$ such that
$$Ent_{\mu}(\mu_t)\leq (1-t)Ent_{\mu}(\rho_0)+tEnt_{\mu}(\rho_1)-\frac{Kt(1-t)}{2}W_{2,\infty}^2(\rho_0,\rho_1)~~~~\forall t
\in [0,1].$$
\end{thm}

For the infinite norm, we will see that $K$ equals to $1$, and for the Sobolev type norm $\|.\|_{k,\gamma}$ that $K$ equals to $1/C_{k,\gamma}^2$.
Precise that Lott-Villani introduced in \cite{Lott} a stronger notion of weak $K-$convexity of relative entropy.
Indeed in its definition, it is required that the property above holds for \underline{all} (constant speed) geodesics.
So in many cases weak $K-$convexity and $K-$convexity coincide, provided that there is unicity of geodesic between two given measures
(this is the case for non branching spaces). An example it fails when optimal coupling is not unique. \\

Now let us briefly summarize the following sections.

Throughout section 2, we resolve the Monge Problem in the Wiener space with the 
cost $\|.\|_{k,\gamma}^p$, in different ways according to parameter $p$. 
Among of these, there is a direct method:
in general settings, the support of an optimal transference plan $\Pi$ (between two probability measures and relatively
to a cost function $c$) is included in $c-$subdifferential of a $c-$convex function (called \emph{potential
of Kantorovich}) $\phi$. It leads to the following system:
\begin{eqnarray*}
\left\{ \begin{array}{llll}
				\phi^c(y)-\phi(x) & = & c(x,y)& \Pi-\textrm{almost everywhere}\\
				\phi^c(y)-\phi(x) &\leq & c(x,y)& \textrm{everywhere}\end{array} \right.
\end{eqnarray*}
And this system can be solved directly when the cost $c$ and the potential $\phi$ 
are differentiable, so long as $\nabla_xc(x,.)$ is injective, 
as it is explained in Villani's book \cite{VIL}. This is the case when $p > 1$. But this method fails
when $p = 1$. In the latter case we head for a recent paper of Bianchini and Cavalletti \cite{BIANCAV}
where the authors resolve Monge Problem in non branching geodesic metric spaces. It turns out that
Wiener space endowed with the norm $\|.\|_{k,\gamma}$ is a such space. 
Simply we will verify suitable conditions.\\

Section 3 is devoted to establish the
$K-$convexity along geodesic of relative entropy (w.r.t. Wiener measure) on the Wiener space
endowed with the infinite norm ($K=1$), then with the Sobolev type norm $\|.\|_{k,\gamma}$ ($K=1/C_{k,\gamma}^2$). 
This time we will process by finite dimensional
approximation as Fang, Shao and Sturm in \cite{FANG}, who have treated the case of the norm $|.|_H$. 
This part requires Wasserstein
distance, which is defined below. Our main contribution consists in establishing results without applying
powerfull tools like \emph{Gromov-Hausdorff} convergence (see \cite{Lott}) or $\mathbb{D}-$convergence
introduced by Sturm in \cite{STURM}. In the language of latter authors, we can say that $(W,|.|_\infty)$
is a $CD(1,\infty)$ space (satisfies such \emph{curvature-dimension condition}) and $(W,\|.\|_{k,\gamma})$
is a $CD(1/C_{k,\gamma}^2,\infty)$ space. Such conditions imply a lot of important results over Wiener spaces. As consequences
over spaces $(W,|.|_\infty),~(W,\|.\|_{k,\gamma})$,
we can quote Brunn-Minkowski, Bishop-Gromov or also Log-Sobolev inequalities (see \cite{USERG}).

\subsection{Settings}

Our ambiant space will be the classical Wiener space $(W,H,\mu)=(\mathcal{C}_0([0,1],\mathbb{R}),H(\mathbb{R}),\mu)$ where
$H(\mathbb{R}):=\{h : [0,1]\longrightarrow \mathbb{R};~h(t)=\int_0^t\dot{h}(s)ds~\mathrm{and}~\dot{h}\in L^2([0,1])\}$ and $\mu$ is the Wiener measure. 
$H(\mathbb{R})$ is a Hilbert space with the inner product:
$$(h,g)_H := \int_{[0,1]}\dot{h}(t)\dot{g}(t)dt.$$
Important facts are that $H$ is dense in $W$ with respect to the uniform norm, and moreover $\mu(H)=0$.

Given two Borel measures $\rho_0$ and $\rho_1$ on $W$, let us state a condition of preserving mass, necessary
for all our discussion:
\begin{equation}\label{PreservingMass}
	\int_W d\rho_0 = \int_W d\rho_1<+\infty.
\end{equation}
In particular it is satisfied when $\rho_0$ and $\rho_1$ are both probability measures. When (\ref{PreservingMass}) is satisfied,
we will consider many times the following Monge-Kantorovich Problem
\begin{equation}\label{MongeKantorovichP}
	\inf_{\Pi\in \Gamma(\rho_0,\rho_1)}\mathcal{I}_p(\Pi)=\inf_{\Pi \in \Gamma(\rho_0,\rho_1)}\int_{W\times W}\|x-y\|_{k,\gamma}^pd\Pi(x,y),
\end{equation}
which is a relaxed problem of Monge Problem (\ref{MongeP}) stated above. At the moment of the cost $c$ is continuous and the ambiant
space is Polish, there is always an existing minimizer for the Monge-Kantorovich Problem (\ref{MongeKantorovichP}) (see for example \cite{AMB}). 
A such minimizer will be called an \emph{optimal coupling} or \emph{optimal transference plan}
between $\rho_0$ and $\rho_1$.\\

We endow the space $W$ with the following norm:
$$\|w\|_{k,\gamma}:=\left(\int_0^1\int_0^1\frac{(w(t)-w(s))^{2k}}{|t-s|^{1+2k\gamma}}dtds\right)^{1/2k}~~~~\mathrm{for}~w\in W,$$
such that $0<\gamma <1/2$, and $2<1+2k\gamma < k$ where $k$ is an integer. In fact this is a pseudo-norm over $W$ since it
can take infinite value. For this reason, we consider $\hat{W}:=\{w\in W;~\|w\|_{k,\gamma}<\infty\}$.
It is well known that $\mu(\hat{W})=1$ hence for a sake of notation, we will write $W=\hat{W}$. $(W,\|.\|_{k,\gamma}$ is
a separable Banach space, and all measures considered in the sequel will be Borel with respect to the topology induced
by the norm $\|.\|_{k,\gamma}$. $H$ is still dense in $(W,\|.\|_{k,\gamma})$.
We can write
\begin{eqnarray*}
	\|w\|_{k,\gamma}^{2k} &\leq & \int_0^1\int_0^1\left(\int_0^1|\dot{x}(\xi)|1\,\,_{s<t}d\xi\right)^2\frac{|x(t)-x(s)|^{2(k-1)}}
		{|t-s|^{1+2k\gamma}}dtds \\
		&\leq & \left(\int_0^1|\dot{x}(\xi)|^2d\xi\right)\int_0^1\int_0^1|t-s|\frac{|x(t)-x(s)|^{2(k-1)}}
		{|t-s|^{1+2k\gamma}}dtds \\
		&\leq & C_{k,\gamma}^{2k}|x|_H^{2k},
\end{eqnarray*}
where $C_{k,\gamma}:=\left(\int_0^1\int_0^1 |t-s|^{k-1-2k\gamma}dtds\right)^{1/2k}$.\\

Before to continue, let us explain how we can decompose the classical Wiener space in finite dimensional spaces.
Consider the projections $\pi_n : W\longrightarrow W$ defined as
$$\pi_n(x)(t):=x\left(\frac{k}{2^n}\right)+2^n\left(t-\frac{k}{2^n}\right)\left[x\left(\frac{k+1}{2^n}\right)-x\left(\frac{k}{2^n}\right)\right],~~ \mathrm{if}~t\in \left[\frac{k}{2^n},
	\frac{k+1}{2^n}\right].$$
At a continuous path, $\pi_n$ associates its affine part.
Denote by $V_n:=\pi_n(W)$. We know that $V_n$ is spanned by the Haar functions and therefore has dimension $2^n$ 
and then will be identified with $\mathbb{R}^N$.

Denote by $V_n^{\bot}$ the subset of $W$ which is image of $W$ by
the map $I_W-\pi_n$. Since $\pi_n$ is a projection, one can write $W=V_n \oplus V_n^{\bot}$.

An important fact is the following: the image measure of the Wiener measure $\mu$ by $\pi_n$ is
the standard Gaussian measure on $V_n$. In other words $(\pi_n)_\#\mu = \gamma_d$, where the 
indice $d$ is the dimension of $V_n$	i.e $2^n = d$. For a sake
	of simplicity, we will denote $(\pi_n)_\#\mu = \gamma_n$ instead of $\gamma_d$.\\

Recall that the \emph{Sobolev space} over the Wiener space $\mathbb{D}_1^p(W)$ is the set of all $\mu-$measurable function
$F\in L^p(W,\mu)$ such that there exists $\nabla F \in L^p(W,H)$ where we have $\mu-$almost surely:
$$(\nabla F(w),h)_H=D_hF(w):=\lim_{\epsilon \rightarrow 0}\frac{F(w+\epsilon h)-F(w)}{\epsilon}~~\mathrm{in}~L^p(W,\mu)~~
\forall h \in H.$$
We will denote 
$$\mathbb{D}_1^\infty(W):=\bigcap_{p>0}\mathbb{D}_1^p(W)~~~~\mathrm{and}~~~~\mathbb{E}(F):=\mathbb{E}_\mu(F).$$

In section 3 we will need $2-$Wasserstein distance between two probability 
measures $\rho_0$ and $\rho_1$ on a measurable space $W$, defined as
\begin{equation}\label{wassersteinDef}
	W_{2,c}^2(\rho_0,\rho_1):=\inf_{\Pi \in \Gamma(\rho_0,\rho_1)}\int c(x,y)^2d\Pi(x,y)
\end{equation}
where $\Gamma(\rho_0,\rho_1)$ is the set of coupling measures between $\rho_0$ and $\rho_1$ i.e the measures
on $W\times W$ with marginals $\rho_0$ and $\rho_1$. 
Moreover we will denote by $\Gamma_0^{c}(\rho_0,\rho_1)$ the set of optimal couplings for $c^2$ or equivalently
the set of couplings which realize minimum on $2-$Wasserstein distance induced by the cost $c$.
A sufficient condition (but not natural)
for that the Wasserstein distance is that $\rho_0$ and $\rho_1$ have second finite moments. To justify terminology,
notice that $W_{2,c}$ is well a distance so long as $c$ is a distance on $W$ and in this case the space
$$\mathcal{P}_2(W):=\left\{\rho ~\mathrm{probability~measure~on~}W;~\int_Wc(x,x_0)^2d\rho(x)<\infty~\mathrm{for~some~}x_0\in W\right\}$$
endowed with $W_{2,c}$ is a metric space.\\

In all this paper, $Pr_i$ (with $1\leq i \leq N$) stand for the projections onto the $i-$th component:
$$Pr_i:X_1\times X_2\times \dots X_N\longrightarrow X_i,$$
where $N\in \mathbb{N}$ depends on the context.

\section{Monge Problem on Wiener space with $\|.\|_{k,\gamma}^p$}

We can now focus on the Monge Problem. 
Notice that $(W,\|.\|_{k,\gamma})$ is a Polish space and in addition
the cost $c(x,y)=\|x-y\|_{k,\gamma}^p$ being continuous (for all $p \geq 1$), 
there is always an existing measure $\Pi \in \Gamma(\rho_0,\rho_1)$ 
which attains the minimum in the Monge-Kantorovich Problem (\ref{MongeKantorovichP}).

For $\rho_0, ~\rho_1$ two measures on $W$ satisfying (\ref{PreservingMass}), recall that the Monge problem between $\rho_0$ and
$\rho_1$ consists of finding an \emph{optimal map} $T:W\longrightarrow W$ which pushes $\rho_0$ forwards to
$\rho_1$ and minimizes the quantity
$$\int_W c(x,G(x))d\rho_0(x),$$
among all push-forward maps $G$ (i.e. $G_\#\rho_0(E):=\rho_0(G^{-1}(E))=\rho_1(E)$ for all Borel $E$ subset of $W$). 

The usual strategy is to use characterization of optimal coupling, with the help of Kantorovich Potentials.
Hence we need the concept of $c-$convexity.

\begin{defin}
Let $\varphi : W \longrightarrow \mathbb{R}$. We say that $\varphi$ is \emph{$c-$convex} if
$$\varphi(x)=\sup_{y\in W}\left(\varphi^c(y)-c(x,y)\right)~~~~\forall x\in W.$$
where $\varphi^c$, called \emph{$c-$transform} of $\varphi$, is defined as:
$$\varphi^c(y)=\inf_{x\in W}\left( \varphi(x)+c(x,y)\right)~~~~\forall y\in W.$$
\end{defin}

Notice that our cost $c$ does not take an infinite value, so $c-$convex functions are real-valued well.

It is well known that any optimal coupling is $c-$cyclically monotone i.e. its \emph{support} 
(defined as the smaller closed subset of $W\times W$ having full $\Pi-$measure) is
$c-$cylically monotone, namely: for all $N\in \mathbb{N}$ and $(x_1,y_1),\dots,(x_N,y_N)\in Supp(\Pi)$ we have
$$\sum_{i=1}^Nc(x_i,y_i)\leq \sum_{i=1}^Nc(x_i,y_{i+1}),$$
with $y_{N+1}:=y_1$.

Rather there is equivalence (in our case) between optimality and $c-$cyclical monotonicity (see e.g.
\cite{USERG}). In particular if $Supp(\Pi)$
is $c-$cyclically monotone then any coupling $\tilde{\Pi}$ such that $Supp(\tilde{\Pi})\subset Supp(\Pi)$ is also $c-$cyclically
monotone.\\

Let us precise what are the difficulties we have met.

\begin{itemize}	
	\item When $p=1$, $c(x,y)=\|x-y\|_{k,\gamma}$ is the considered norm on $W$. Hence if a map $\varphi$ is 
$c-$convex then it is also $1-$Lipschitz, hence $H-$Lipschitz. Indeed
$$|\varphi(x+h)-\varphi(x)|\leq \|h\|_{k,\gamma}\leq C_{k,\gamma}|h|_H ~~~~\forall h\in H~~\forall x\in W.$$
In this case, we have a version of Rademacher theorem (proved in \cite{ENCHSTRO} and recalled  in appendice
for a sake of completeness) on Wiener space, for such $H-$Lipschitz functions.
With the remark above, it leads to $c-$convex functions on the Wiener space are almost surely differentiable.
But the difficulty in this case is that the cost, being a norm, is \emph{not strictly convex}, so we do not have
injectivity of $y\longmapsto \nabla_xc(x,y)$, therefore our method is not available.
Nevertheless we shall use an other method to solve Monge Problem. Indeed thanks to the part $2.$ of 
Lemma (\ref{PropertiesN}), $(W,\|.\|_{k,\gamma})$ is a geodesic non branching space (i.e.
geodesics cannot bifurcate). So we can apply method,
detailed in \cite{BIANCAV}.\\

	\item When $p>1$, $c(x,y)=\|x-y\|_{k,\gamma}^p$ becomes strictly convex, we get this time the injectivity
of $\nabla_xc(x,.)$. But we lose the $H-$Lipschitz property of
$c-$convex functions. Indeed if $\varphi$ is such function we can write
\begin{eqnarray*}
	|\varphi(x)-\varphi(y)|&\leq& |\|x-\xi\|_{k,\gamma}^p-\|y-\xi\|_{k,\gamma}^p|\\
		&\leq & \|x-y\|_{k,\gamma}M_{\xi},
\end{eqnarray*}
where the latter constant $M_{\xi}$ depends on $\xi$ and cannot be bounded.
However we will see that in this case $c-$convex functions (hence potentials) are
\emph{locally} $H-$Lipschitz. Since differentiability is a local property, we should be able to apply 
Rademacher theorem again.
\end{itemize}

Before we continue, let us set properties of the considered norm, which will be useful for the sequel.

\subsection{Properties of $\|.\|_{k,\gamma}$}

We give two ingredients that will be essential for the sequel.

\begin{lemme}\label{PropertiesN}
If we denote by $\tilde{F}:W\longrightarrow \mathbb{R}_+$ the map $\tilde{F}(w)=\|w\|_{k,\gamma}$, then we have the following properties:
\begin{enumerate}
	\item $\tilde{F}$ admits a gradient $\nabla \tilde{F}(w)$ belonging to  $W^\star$ for all $w\in W\backslash \{0\}$, where
	$W^\star$ is the dual of $W$. Moreover $\tilde{F}^p$ is everywhere differentiable for all $p>1$.
	\item $\tilde{F}$ is a norm such that its unit ball is strictly convex.
\end{enumerate}
\end{lemme}

The first part of the proof is inspired from \cite{FANG0}. 

\begin{demo}
\noindent $1.$ First we show the property for $F:=\tilde{F}^{2k}$. 
Take $h\in W$, we can write for $w \in W$ and $\epsilon >0$,
$$F(w+\epsilon h)=\int_0^1\int_0^1\frac{\left((w(t)-w(s))+\epsilon(h(t)-h(s))\right)^{2k}}{|t-s|^{1+2k\gamma}}dtds.$$
And taking the derivative at $\epsilon = 0$, it is clear that $\lim_{\epsilon \rightarrow 0}\frac{F(w+\epsilon h)-F(w)}{\epsilon}$
exists and moreover:
\begin{eqnarray*}
	|D_hF(w)|&\leq &2k\int_0^1\int_0^1\frac{|w(t)-w(s)|^{2k-1}}{|t-s|^{1+2k\gamma}}|h(t)-h(s)|dtds\\
		&\leq &2k\int_{[0,1]^2}\frac{|w(t)-w(s)|^{2k-1}}{|t-s|^{(1+2k\gamma)(2k-1)/(2k)}}
		\frac{|h(t)-h(s)|}{|t-s|^{(1+2k\gamma)/(2k)}}dtds
\end{eqnarray*}
and now applying Hölder's inequality, we get
\begin{eqnarray*}
	|D_hF(w)|&\leq &2k\left(\int_{[0,1]^2}\frac{|w(t)-w(s)|^{2k}}{|t-s|^{1+2k\gamma}}dtds\right)^{(2k-1)/(2k)}
		\left(\int_{[0,1]^2}\frac{|h(t)-h(s)|^{2k}}{|t-s|^{1+2k\gamma}}dtds\right)^{1/(2k)}\\
		&=&2k\|w\|_{k,\gamma}^{2k-1}.\|h\|_{k,\gamma}.
\end{eqnarray*}
Hence $h\longmapsto D_hF(w)$ is a bounded operator on $W$ for all $w\in W$. It leads to existence of a
gradient $\nabla F(w)$ which belong to the dual space $W^\star \subset H^\star=H$ (by (\ref{inqNorms})). Since $\tilde{F}=F^{1/(2k)}$,
its gradient satisfies $\nabla \tilde{F}(w)=F^{1/(2k)-1}(w)\nabla F(w)$ and in particular
$w$ must not be equal to zero.

Since $\tilde{F}$ is differentiable out of $\{0\}$ it is just a remark to see that for any $p>1$,
$\tilde{F}^p$ is differentiable everywhere over $(W,\|.\|_{k,\gamma})$.

\noindent $2.$ This proof is the same as the proof of Minkowski's inequality.
Indeed if $w_1,~w_2\in W$ and $\eta \in (0,1)$ then we have
\begin{eqnarray*}
	\|(1-\eta )w_1+\eta w_2\|_{k,\gamma}^{2k}
	&=&\int_{[0,1]^2}\frac{|(1-\eta )(w_1(t)-w_1(s))+\eta(w_2(t)-w_2(s))|
	^{2k}}{|t-s|^{1+2k\gamma}}dtds\\
	& =& \int_{[0,1]^2}|(1-\eta )(w_1(t)-w_1(s))+\eta(w_2(t)-w_2(s))|\\
	&\times &\frac{|(1-\eta )(w_1(t)-w_1(s))+\eta(w_2(t)-w_2(s))|
	^{2k-1}}{|t-s|^{1+2k\gamma}}dtds\\
	&\leq & \int_{[0,1]^2}\frac{(1-\eta )|w_1(t)-w_1(s)|}{|t-s|^{(1+2k\gamma)/(2k)}}\frac{
	|(1-\eta )(w_1(t)-w_1(s))+\eta(w_2(t)-w_2(s))|^{2k-1}}{|t-s|^{(1+2k\gamma-\frac{1}{2k}-\gamma)}}dtds\\
	&+ &	\int_{[0,1]^2}\frac{\eta |w_2(t)-w_2(s)|}{|t-s|^{(1+2k\gamma)/(2k)}}\frac{
	|(1-\eta )(w_1(t)-w_1(s))+\eta(w_2(t)-w_2(s))|^{2k-1}}{|t-s|^{(1+2k\gamma-\frac{1}{2k}-\gamma)}}dtds\\
	&\leq & \left((1-\eta)\|w_1\|_{k,\gamma}+	\eta \|w_2\|_{k,\gamma} \right)\left(\|(1-\eta)w_1+\eta w_2\|_{k,\gamma}^{2k}
	\right)^{1-1/2k}.
\end{eqnarray*}
The two inequalities above are respectively triangle inequality and Holder's inequality, and are in fact 
equality if and only if
$w_1$ and $w_2$ are almost everywhere colinear. This leads to the strict convexity of our norm.
\end{demo}

Conditions on parameters $(p,\gamma)$ are sufficient to have $F \in \mathbb{D}_1^\infty(W)$, as it is shown
in \cite{FANG0}. Notice that the (Gateaux) differentiability of $\tilde{F}$ exists in the direction of $W$, hence in
particular in the direction of $H$.

\subsection{The case $p>1$.}

Throughout this subsection, the cost is $c(x,y)=\|x-y\|_{k,\gamma}^p$, with $p>1$.

We follow Fathi and Figalli in \cite{FAFI} to get around the fact that $c-$convex functions are not
$1-$Lipschitz with respect to $\|.\|_{k,\gamma}$, but nevertheless are locally Lipschitz with restriction to
suitable subsets.
The key argument is that the $\sup$ of a family of uniformly $\|.\|_{k,\gamma}-$Lipschitz functions, 
is also $\|.\|_{k,\gamma}-$Lipschitz.

\begin{thm}\label{MongeSolveKGAMMA}
Let $\rho_0$ and $\rho_1$ be two measures on $W$ satisfying (\ref{PreservingMass}) and such that the first one is absolutely
continuous with respect to the Wiener measure $\mu$.
Assume $\mathcal{I}(\Pi)$ is finite for some $\Pi\in \Gamma(\rho_0,\rho_1)$. 

Then there exists a unique optimal coupling between $\rho_0$ and $\rho_1$ relatively to the cost $c$.
Moreoever it is concentrated on a graph of some Borel map $T:W\longrightarrow W$
unique up to a set of zero measure for $\mu$. 
\end{thm}

\begin{demo}
Let $\Pi \in \Gamma(\rho_0,\rho_1)$ be an optimal coupling for $c$. We shall show that $\Pi$ is concentrated on
a graph of some Borel map. 
It is well known  (see e.g. \cite{VIL}) that under condition $\mathcal{I}(\Pi)$ is finite, since $Supp(\Pi)$ is $c-$cyclically monotone,
there is a $c-$convex map $\varphi :W \longrightarrow \mathbb{R}$ (called Kantorovich potential) such that
$$	\varphi^c(y)-\varphi(x)=\|x-y\|_{k,\gamma}^p~~~~\Pi-\mathrm{a.s.}$$
Moreover from the definition of $c-$convexity, we also have
\begin{equation}\label{subdifferential}
	\varphi^c(y)-\varphi(x)\leq \|x-y\|_{k,\gamma}^p~~~~\forall (x,y)\in W\times W.
\end{equation}
Since $\varphi^c$ is finite everywhere, if we consider subsets $W_n:=\{\varphi^c \leq n\}$ for $n\in \mathbb{N}$ then
$$W_n\subset W_{n+1}~~\mathrm{and}~~\bigcup_{n\in \mathbb{N}} W_n = W.$$
Our cost $c(.,y)=\|.-y\|_{k,\gamma}^p$ is locally $\|.\|_{k,\gamma}-$Lispchitz locally uniformly in $y$.
Hence for each $y\in W$ there exists a neighborhood $E_y$ of $y$ such that $(\|.-z\|_{k,\gamma}^p)_{z\in E_y}$
is a uniform family of locally $\|.\|_{k,\gamma}-$Lipschitz functions.
Moreover $W$ being separable, we can find a sequence $(y_l)_{l\in \mathbb{N}}$ of elements of $W$ such that
$$\bigcup_{l\in \mathbb{N}} E_{y_l} = W.$$
Now consider increasing subsets of $W$:
$$V_n:=W_n \bigcap (\bigcup_{l=1}^n E_{y_l}).$$
We can define maps approximating $\varphi$ as follow
\begin{eqnarray*}
	\varphi_n : W & \longrightarrow & W\\
	x & \longmapsto & \sup_{y\in V_n} \left(\varphi^c(y)-\|x-y\|_{k,\gamma}^p\right).
\end{eqnarray*}
Notice that 
$$\varphi_n(x)=\max_{l=1,\dots,n} \sup_{y\in W_n\cap E_{y_l}}\left(\varphi^c(y)-\|x-y\|_{k,\gamma}^p\right).$$
But since $\varphi \leq n$ on $W_n$ and $-\|.\|_{k,\gamma}^p$ is bounded from above, $\varphi_n$ is also bounded from
above.
Therefore the sequel $(\varphi^c(y)-\|.-y\|_{k,\gamma}^p)_{y\in W_n \cap E_{y_l}}$ is uniformly locally 
$\|.\|_{k,\gamma}-$Lipschitz and bounded
from above. Finally Proposition \ref{SupConvex} in Appendices shows that $\varphi_n$ being a maximum of uniformly
locally $\|.\|_{k,\gamma}-$Lipschitz
functions, is also locally $\|.\|_{k,\gamma}-$Lispchitz. We can extend $\varphi_n$ to a $\|.\|_{k,\gamma}-$
Lipschitz function everywhere on $W$ still denoted by $\varphi_n$. By (\ref{inqNorms}), we get
$$|\varphi_n(w+h)-\varphi_n(w)|\leq C\|h\|_{k,\gamma}\leq 2C|h|_H~~\forall w\in W,~\forall h \in H.$$
Namely, $\varphi_n$ is a $H-$Lipschitz function. Thanks to Rademacher
theorem, there exists a Borel subset $F_n$ of $W$ with plain $\mu-$(hence $\rho_0-$)measure such that 
for all $x\in F_n$, $\varphi_n$ is differentiable at $x$. Then for each $x \in F:=\cap_n F_n$ (which has also plain
$\rho_0-$measure),
each $\varphi_n$ is differentiable at $x$.

By increasing of $(V_n)_n$, it is clear that $\varphi_n \leq \varphi_{n+1} \leq \varphi$ everywhere on $W$.
Moreover with same argument as in \cite{FAFI}, if $P_n:=Pr_1\left(Supp(\Pi)\cap(W\times V_n)\right)$ then
$\varphi_{|P_n}=\varphi_{n|P_n}=\varphi_{l|P_n}$ for all $l \geq n$ and all $n \in \mathbb{N}$.
Fix $x\in P_n \cap F$. By definition of $P_n$ it exists $y_x\in V_n$ such as
\begin{eqnarray*}
	\varphi^c(y_x)-\varphi_n(x)&=&\|x-y_x\|_{k,\gamma}^p,\\
	\mathrm{i.e.}~\varphi^c(y_x)-\varphi (x)&=&\|x-y_x\|_{k,\gamma}^p.
\end{eqnarray*}
Subtracting (\ref{subdifferential}) with $(x',y_x)$ to the previous equality, we get for all $x'\in W$ and $h\in H$:
$$\varphi(x)-\varphi(x') \geq \|x-y_x\|_{k,\gamma}^p-\|x'-y_x\|_{k,\gamma}^p.$$
Taking $x'=x+\epsilon h$ with $\epsilon>0$, $h\in H$, dividing by $\epsilon$ and taking the limit when $\epsilon$ tends to $0$,
we get (by linearity in $h$):
\begin{equation}\label{MongeEGrad}
	\nabla \varphi (x) +\nabla_xc(x,y_x)=0.
\end{equation}
Indeed $c(.,y_x)$ is differentiable at $x$ thanks to Proposition \ref{PropertiesN}. 
The strict convexity of $c(x,y)=\|x-y\|_{k,\gamma}^p$ yields $\nabla_xc(x,.)$ is injective and 
(\ref{MongeEGrad}) gives:
$$y_x=(\nabla_xc(x,.))^{-1}(-\nabla \varphi(x))=:T(x),$$
where $(\nabla_xc(x,.))^{-1}$ is the inverse of the map $y\longmapsto \nabla_xc(x,y)$.
Notice here that $T$ is uniquely determined.
We deduce that $Supp(\Pi)\cap (W\times V_n)$ is the graph of the map $T$ over $P_n\cap F$ for all $n \in \mathbb{N}$. 
But $(P_n)_n$ and $(V_n)_n$ are increasing and such that $\bigcup_n V_n = W$.
Therefore $Supp(\Pi)$ is a graph over $Pr_1(Supp(\Pi))\cap F$ with $Pr_1(Supp(\Pi))=\bigcup_n P_n$.

We can extend $T$ onto a measurable map over $W$ as it is explained in \cite{FAFI}.
We obtain $Supp(\Pi)$ is included in the graph of a measurable map $T$, unique up to a set of $\rho_0-$measure. 
In other words $\Pi=(id\times T)_\#\rho_0$.

We have proved that any optimal coupling is carried by a graph of some map. So if $\Pi_1$, 
$\Pi_2\in \Gamma(\rho_0,\rho_1)$ are optimal for $\|.\|_{k,\gamma}$ then any convex combination of $\Pi_1$ and $\Pi_2$
is also optimal. Take $\Pi:=\frac{1}{2}(\Pi_1+\Pi_2)$ be an optimal coupling between $\rho_0$ and $\rho_1$: there
exists some measurable map $T$ such that $\Pi=(Id\times T)_\#\rho_0$. Let $f$ be the density of $\Pi_1$ with respect
to $\Pi$. Then for any continuous bounded functions $\varphi$ we have
\begin{eqnarray*}
	\int_W \varphi(x)d\rho_0(x) &=&	\int_{W\times W}\varphi(x)d\Pi_1(x,y) \\
		&=& \int_{W\times W}\varphi(x)f(x,y)d\Pi(x,y)\\
		&=&\int_W \varphi(x)f(x,T(x))d\rho_0(x).
\end{eqnarray*}
This yields $f(x,T(x))=1$ $\rho_0-$a.e., hence $f=1$ $\Pi-$a.e. It leads to $\Pi=\Pi_1$ and finally $\Pi_2=\Pi_1=(Id\times T)_\#\rho_0$.
\end{demo}

\subsection{The case $p=1$.}

Throughout this subsection, the cost is $c(x,y)=\|x-y\|_{k,\gamma}$.

We follow the method of \cite{CAVA} developed by Bianchini and Cavalletti in \cite{BIANCAV}. For this we will need
notions of \emph{transport rays}, \emph{transport sets} and \emph{endpoints}.
By strict convexity of our norm $\|.\|_{k,\gamma}$, $(W,\|.\|_{k,\gamma})$ is a geodesic non branching space.

Let us recap briefly the different steps of their method:
\begin{enumerate}
	\item	reduce the initial Monge-Kantorovich Problem to the (one-dimensional) Monge-Kantorovich 
	Problem along distinct geodesics (this is possible
		since the space is non-branching).
	\item	verify that the {conditional measures} provided by 
	disintegration of both measures $\rho_0$ and $\rho_1$ on each geodesic are without atoms (this is possible 
		thanks to properties of Gaussian measure) in order to apply result of one-dimensional Monge Problem: on each geodesic
		there exists a transport map. 
	\item piece obtained maps together to get a transport map for the initial Monge Problem (this is possible by a general
		selection theorem).
\end{enumerate}

We present only sketches of results, which are very similar with the Cameron-Martin norm $|.|_H$.
For all details, consult \cite{CAVA} and \cite{BIANCAV}. In our case, the cost $\|.\|_{k,\gamma}$
is smooth enough to guarantee the existence of a Potential $\varphi$ such that for all optimal coupling $\Pi$ we have:
\begin{eqnarray*}
	Supp(\Pi) \subset \Gamma:=\{(x,y)\in W\times W;~\varphi^c(y)-\varphi(x)=\|x-y\|_{k,\gamma}\}.
\end{eqnarray*}

\begin{defin}
 The set of \emph{oriented transport rays} is defined as:
 $$G:=\left\{(x,y)\in W\times W; ~\exists (w,z)\in \Gamma;~\|w-x\|_{k,\gamma}+\|x-y\|_{k,\gamma}+\|y-z\|_{k,\gamma}=
 	\|w-z\|_{k,\gamma}\right\}.$$
 We denote by $G(x)$ the \emph{outgoing transport rays} from $x\in W$ and by $G^{-1}(x)$ the \emph{incoming transport
 rays} in $x\in W$. Finally define the set of \emph{transport rays} as
 $$R:=G\cup G^{-1}.$$
\end{defin}

\begin{defin}
	The \emph{transport set (with endpoints)} is defined as
	$$\mathcal{T}_e:=\{x \in W;~R(x)\neq \{x\}\}.$$
	The \emph{transport set (without endpoints)} is defined as
	$$\mathcal{T}:=\{x\in W;~G(x)\neq \{x\},~G^{-1}(x)\neq \{x\}\}.$$
\end{defin}
	
Here the important fact is the following: the space $(W,\|.\|_{k,\gamma})$ is non branching hence
the set $R(x)$ is a single geodesic for each $x \in \mathcal{T}$. It leads to the fact that 
$R\cap \mathcal{T}\times \mathcal{T}$ is an equivalence relation on $\mathcal{T}$.
Equivalence classes form a partition in $\mathcal{T}$ and therefore
we can apply a theorem of disintegration of measures, as it is explained in \cite{CAVA}, provided
there is a $\rho_0$-measurable \emph{cross section} 
$f:\mathcal{T}\longrightarrow \mathcal{T}$ for the ray equivalence relation $R$ (see \cite{CAVA} for the terminology).

This is possible to follow the part 4. of \cite{BIANCAV} to prove the existence of such cross section $f$,
since geodesics of $(W,\|.\|_{k,\gamma})$ are continuous and locally compact.

\begin{defin}
	The multivalued \emph{endpoint graphs} are defined as
\begin{eqnarray*}
	a&:=&\{(x,y)\in G^{-1}; ~G^{-1}(y)\backslash \{y\}= \emptyset\}\\
	b&:=&\{(x,y)\in G;~G(y)\backslash \{y\}=\emptyset\}.
\end{eqnarray*}
	Then $Pr_2(a)$ and $Pr_2(b)$ are respectively called \emph{initial points} and \emph{final points}.
\end{defin}	

Since from a point of $a(\mathcal{T})$, many geodesics can start, proposition \ref{InitialFinalNul}
is necessary to see in fact that such set has $\rho_0-$measure null. It is the same for $b(\mathcal{T})$.
Hence the transport rays provide a partition of $\mathcal{T}$ up to a $\mu-$negligible set, and
we disintegrate the measure $\mu$ on $\mathcal{T}_e$ w.r.t. to this partition.
Proposition \ref{AbsContP1} says that conditional measures have no atoms, hence we get an optimal
map for the induced Monge Problem.\\

From now, let us consider an \emph{optimal} transference plan $\Pi_0$ between $\rho_0$ and $\rho_1$ relatively
to the cost $c(x,y)=\|x-y\|_{k,\gamma}$. Denote by $\rho_0^n:=(\pi_n)_\#\rho_0$ and $\rho_1^n:=(\pi_n)_\#\rho_1$. 
This latter cost is a norm over $V_n:=(\pi_n)_\#W$ which inherits strict convexity
and differentiability.
So by the \cite{CHAMP}, the Monge Problem in these settings admits at least one solution, say $T_n$
and $\Pi_0^n:=(id\times T_n)_\#\rho_0^n$ is the unique optimal transference plan between $\rho_0^n$ and $\rho_1^n$.
In other words, $\Pi_0^n$ is concentrated on some Borel set $\Gamma_n\subset Graph(T_n)$.

\vskip 2mm

First we have two technical propositions.

\begin{prop}\label{finitedimest}
Assume that there exists $M>0$ such that densities of $\rho_0^n$ and $\rho_1^n$ are bounded by $M$
$\gamma_n-$almost everywhere. Then the following estimate holds true for all compact subset $A\subset W$:
$$\gamma_n(T_{n,t}(A))\geq \frac{1}{M}\rho_0^n(A)~~~~~~\forall t\in [0,1],$$
where $T_{n,t}:=(1-t)Id+tT_n$.
\end{prop}
	
The proof is quite the same as in \cite{CAVA}. Only difference is to consider Monge maps for the cost
induced by $\|.\|_{k,\gamma}^p$ with ($p>1$), instead of $|.|^p$. Indeed costs $\|.\|_{k,\gamma}^p$ satisfy
conditions of Theorem 6.2.7 in \cite{AMB}, so well that optimal maps $T_n^p$ are approximately differentiable.

Come back to the Wiener space. We have the following approximation result:

\begin{prop}\label{finitetoinfiniteest}
	Assume that there exists $M>0$ such that for all compact subset $A\subset W$,
	the following holds true:
	$$\gamma_n(T_t(\Gamma_n\cap A\times W))\geq M\rho_0^n(A)~~~~\forall n\in \mathbb{N}.$$
	Then for all compact subset $A\subset W$, we have
	$$\mu(T_t(\Gamma\cap A\times W))\geq M\rho_0(A),$$
	where $T_t(x,y):=(1-t)x+ty$.
\end{prop}

The proof (established again in \cite{CAVA})
uses generalities of measure theory, Hausdorff topology. It is true with general measures,
provided the cost is at least lower semi-continuous.

We can head for the solution of Monge Problem.

The first step is to prove that set of initial points has $\rho_0-$measure zero and final points has
$\rho_1-$measure zero. We denote for the sequel by $m:=f_\#\rho_0$ and $m_\mu:=f_\#\mu$ the image measures
by $f$ (being the cross section).

\begin{prop}\label{InitialFinalNul}
	 If $\rho_0$ and $\rho_1$ are absolutely continuous with respect to the Wiener measure $\mu$ then:
	 $$\rho_0(a(\mathcal{T}))=0~~~~~~\mathrm{and}~~~~~~\rho_1(b(\mathcal{T}))=0.$$
\end{prop} 

\begin{demo}
	We can only prove that $\rho_0(a(\mathcal{T}))=0$.
	Let $f_0$ be the density of $\rho_0$ with respect to $\mu$.
	Assume that $\rho_0(a(\mathcal{T}))>0$.
	Let $A\subset a(\mathcal{T})$ and $\delta, M>0$ such that $\rho_0(A)>0$ and for every $x\in A$, $\delta < f_0(x)\leq M$.
	Consider the restriction $\mu_{|\mathcal{T}}$ of $\mu$ with respect to $\mathcal{T}$ and its disintegration: 
	$$\mu_{|\mathcal{T}}=\int_{\mathcal{S}} \mu_y m_\mu (dy),~~~~\mu_y(\mathcal{T})=1~~m_\mu-\mathrm{a.e.}$$
	Now consider the initial point map $a:\mathcal{S}\longrightarrow A$ and the measure image $a_\#m_\mu$.
	We observe that $\rho_{0|A}$ is absolutely continuous with respect to $a_\#m_\mu$. Indeed we have:
	$$\forall B\subset A,~\rho_0(B)>0 \Rightarrow \mu(R(B)\cap \mathcal{T}) >0.$$
	Hence there exists a subset $\hat{A}\subset A$ of positive $a_\#m_\mu-$measure such that the map $h:\hat{A}\longrightarrow
	\mathbb{R}$ defined as:
	$$h(x):=\frac{d\rho_{0|A}}{da_\#m_\mu}(x)$$
	satisfies $h(x)\leq M'$ for some positive constant $M'$.
	
	Now let us introduce
	$$\hat{\mu}(.):=\int_{R(\hat{A})\cap \mathcal{S}} h(a(y))\mu_y(.) dm_\mu (dy).$$
	For $S\subset \mathcal{S}$ we have:
	\begin{eqnarray*}
		\rho_{0|\hat{A}}(\cup_{y\in S}R(y)) &=& \rho_{0|\hat{A}}(a(S))\\
			&=&\int_{a(S)}h(w)d(a_\#m_\mu)(w)\\
			&=&\int_Sh(a(w))m_\mu(dw)= \hat{\mu}(\cup_{y\in S}R(y)).
	\end{eqnarray*}
	This yields that $\mathcal{T}_e$ is still a transport set for the transport problem between
	$\rho_{0|\hat{A}}$ and $\hat{\mu}$. Moreover these two measures have uniformly bounded densities.
	Hence if we project them, we obtain the finite dimensional estimate of Proposition \ref{finitedimest}.
	Then we get the infinite dimensional estimate thanks to Proposition \ref{finitetoinfiniteest}:
	there exists some $C>0$ such that for any compact subset $B$ of $W$, we have:
	$$\mu(T_t(\Gamma \cap B\times W))\geq C\rho_0(B).$$
	This is the case for $A$ considered above. Then we can find a sequence $(t_n)_n$ converging to $0$ such that:
\begin{eqnarray*}
	\mu(T_{t_n}(\Gamma \cap A\times W))&\geq& C\rho_0(A) \geq \delta C\mu(A)\\
	\mathrm{and}~~ A\cap T_{t_n}(\Gamma \cap A\times W)&=&\emptyset.
\end{eqnarray*}
	If we denote by $A^\epsilon:=\{x\in W;~\|A-x\|_{k,\gamma} < \epsilon \}$ then for $t_n \leq \epsilon /M$
	we have:
	$$\mu(A^\epsilon)\geq \mu(A)+\mu(T_{t_n}(\Gamma \cap A\times W)) \geq (1+C\delta)\mu(A).$$
	But $\mu(A)=\lim_{\epsilon \rightarrow 0}\mu(A^\epsilon)$ and we get a contradiction. The result follows.
\end{demo}

Here we state a fundamental point which will allow us to apply result of Monge Problem in one dimension.

\begin{prop}\label{AbsContP1}
	If $\rho_0$ and $\rho_1$ are absolutely continuous with respect to the Wiener measure $\mu$ then 
	for $m-$a.e. $y\in \mathcal{S}$, the conditional probabilities $\rho_{0,y}$ and $\rho_{1,y}$ have no atoms. 
\end{prop}

\begin{demo}
	We can only prove for the first measure.
	Let $f_0$ and $f_1$ be respectively densities of $\rho_0$ and $\rho_1$ with respect to $\mu$.
	Assume that there exists a measurable subset $\hat{S}$ of $\mathcal{S}$ such that $m(\hat{S})>0$ and for every
	$y\in \hat{S}$ there exists $x_y$ such that $\rho_{0,y}(\{x_y\})>0$.
	From now we restrict both $\rho_0$ and $\rho_1$ to $R(\hat{S})$ and denote again with $\rho_0$ and $\rho_1$.
	Let us consider the sets $K_{i,M}:=\{x\in X;~f_i\leq M\}$ for $i=0,1$. For $M$ sufficiently large,
	the conditional probabilities of the disintegration of $\rho_{0|K_{0,M}}$ have atoms.
	
	Define $\rho_{0,y,M}:=\rho_{0,y|K_{0,M}}$, $\rho_{1,y,M}:=\rho_{1,y|K_{1,M}}$ and 
	$$D(N):=\left\{y\in \hat{S}; ~h(y):=\frac{\rho_{0,y,M}(R(y))}{\rho_{1,y,M}(R(y))}\leq N \right\}.$$
	For $N$ sufficiently large, we have $m(D(N))>0$. Hence the map $h: D(N)\longrightarrow \mathbb{R}$
	is well defined and permits to introduce
	$$\hat{\rho_0}:=\int_{D(N)}h(y)\rho_{0,y,M}m(dy),~~~~~~\hat{\rho_1}:=\rho_{1|R(D(N))\cap K_{1,M}}.$$
	These measures have bounded densities with respect to $\mu$ and the set $\hat{\mathcal{T}}:=
	\mathcal{T}\cap G(K_{0,\delta})	\cap G^{-1}(K_{1,\delta})$ 
	is a transport set for the transport problem between $\hat{\rho_0}$ and $\hat{\rho_1}$.
	It follows that $\hat{\mu}:=\mu_{|\hat{\mathcal{T}}}$ satisfies conditions of Proposition \ref{finitetoinfiniteest}.
	
	Suppose $\hat{\rho_0}(Pr_2(A))>0$ where
	$$A:=\bigcup_{y\in \mathcal{S}} \{x\in R(y); ~\hat{\rho_0}_y(\{x\})>0\}$$
	is a Borel set. It is well known that $A$ is a countable union of Borel graphs (Lusin Theorem).
	If we take one of such graph $\hat{A}$ we have $m(Pr_1(\hat{A}))>0$ hence by disintegration
	$\rho_0(Pr_2(\hat{A}))>0$. Applying Proposition \ref{finitetoinfiniteest} we get:
	$$\mu (T_{t}(\Gamma \cap Pr_2(\hat{A})\times W)) \geq M \rho_0(Pr_2(\hat{A})).$$
	But $T_t(\Gamma \cap Pr_2(\hat{A})\times W) \cap (Pr_2(\hat{A})) =\emptyset$ then letting $t\longrightarrow 0$
	we get a contradiction and it follows $\hat{\rho_0}(Pr_2(A))=0$.
	
	In particular the conditional probabilities $\hat{\rho_0}_y$ of the disintegration of $\hat{\rho_0}$
	have no atoms. We get a contradiction since $\hat{\rho_0}_y$ is absolutely continuous with respect to
	$\rho_{0,y|\hat{\mathcal{T}}}$ and the latter measure has atoms.
\end{demo}

\begin{thm}
Let $\rho_0$ and $\rho_1$ be two measures on $W$ satisfying (\ref{PreservingMass}) and such that both are absolutely
continuous with respect to the Wiener measure $\mu$.
Assume $\mathcal{I}(\Pi)$ is finite for some $\Pi\in \Gamma(\rho_0,\rho_1)$. 

Then there exists an optimal coupling between $\rho_0$ and $\rho_1$ which is
concentrated on a graph of some Borel map $T:W\longrightarrow W$.
\end{thm}

\begin{demo}
	We can process as Theorem 8.3 in \cite{CAVA}, putting $\|.\|_{k,\gamma}$ instead of $|.|_H$.
\end{demo}

\section{Convexity of relative entropy on Wiener space}

This part is split in two sections. The first one contains recalls in finite dimensional cases and an
extension of known results for uniform and $\|.\|_{k,\gamma}$ norms on $\mathbb{R}^d$. 
The second one will need this result to show our
purpose on the Wiener space. Throughout all of this section, the cost $c$ will be a distance induced by a norm.
So that it takes sense to consider Wasserstein distance.

We consider the relative entropy with respect to a reference measure $\gamma$ defined as
\begin{eqnarray}
	Ent_{\gamma}(\rho) := \left\{ \begin{array}{ll}
				\int f \log (f) d\gamma & \textrm{if $\rho$ admits $f$ for density w.r.t $\gamma$}\\
				+\infty& \textrm{otherwise}\end{array} \right.
\end{eqnarray}
We need to recall the notion of \emph{geodesics} on the space of probability measures.
$t \in [0,1] \longmapsto \rho_t^{(c)} \in \mathcal{P}_2(W)$ is a \emph{(constant speed) geodesic}, provided
$$W_{2,c}(\rho_t^{(c)},\rho_s^{(c)})=(t-s)W_{2,c}(\rho_0,\rho_1)~~~~\forall 0\leq s \leq t\leq 1.$$
One can obtain a constant speed geodesic by letting $\rho_t^{c}:=((1-t)Pr_1+tPr_2)_{\#}\Pi,$ $\forall t\in [0,1],$
if $\Pi$ is an optimal coupling for $c$. In fact each optimal transference plan involves
a constant speed geodesic (see for instance \cite{USERG}).
Moreover on non branching space, whenever the optimal transference plan is unique, 
there is a \emph{unique} constant speed geodesic between $\rho_0$ and $\rho_1$. For example, Banach space with strictly
convex norm is non branching, while Banach space with non strictly convex norm is branching.

Thanks to this definition, one can consider the notion of \emph{convexity} (and \emph{weak convexity}) 
\emph{along (constant speed) geodesics}.

\begin{defin}
For $\rho_0,~\rho_1$ two probability measures with second finite moments (for a cost $c$), 
we say relative entropy with respect
to a reference measure $\mu$, is \emph{weakly $K-$convex along (constant speed) geodesics on $(\mathcal{P}_2(.),c)$} if
there exists a (constant speed) geodesic
$\rho_t^{(c)}$ induced by an optimal transference plan belonging to $\Gamma_0^c(\rho_0,\rho_1)$ such that
$$Ent_{\mu}(\rho_t^{(c)})\leq (1-t)Ent_{\mu}(\rho_0)+tEnt_{\mu}(\rho_1)-\frac{Kt(1-t)}{2}W_{2,c}^2(\rho_0,\rho_1)~~~~\forall t
\in [0,1].$$
We say that relative entropy is \emph{$K-$convex along (constant speed) geodesics on $(\mathcal{P}_2(.),c)$} 
(\underline{not weakly}) if the latter inequality holds \underline{for all} (constant speed) geodesics $\rho_t^{(c)}$.
\end{defin}

\subsection{Finite dimensional cases}

We want to show the relative entropy (with respect to the standard Gaussian measure $\gamma_d$) is $K-$convex along
(constant speed) geodesics in the finite dimensional space $\mathbb{R}^d$ endowed with the uniform norm and the norm $\|.\|_{k,\gamma}$.
For the first one, we need to consider
the cases where $\mathbb{R}^d$ is endowed with different $p-$norms.
For $p$ positive integer we set:
$$c_p(x,y):=|x-y|_p=\left(\sum_{i=1}^d|x_i|^p\right)^{1/p}.$$
Let us begin to recall $K-$convexity along geodesics for suitable norm on $\mathbb{R}^d$.  

\begin{prop}\label{EntRnNorm}
	Let $\|.\|$ be a strictly convex and differentiable norm on $\mathbb{R}^d\backslash \{0\}$. If $\sqrt{K} \|.\|\leq |.|_2$
	then relative entropy w.r.t. $\gamma_d$ on $(\mathbb{R}^d,\|.\|)$ is $K-$convex 
	along (constant speed) geodesics on $(\mathcal{P}_2(\mathbb{R}^d),W_{2})$, where $W_2$ is the Wasserstein
	distance induced by $\|.\|$.
\end{prop}

By assumption on $\|.\|$, here $(\mathbb{R}^d,\|.\|)$ is a non branching space. 

\begin{demo}
Let $\rho_0$ and $\rho_1$ be two probability measures absolutely continuous with respect to $\gamma_d$ (hence $\mathcal{L}$) 
(otherwise the result is trivial). For $i=0,1$ let $d\rho_0=f_0d\mathcal{L}$ and $d\rho_1=f_1d\mathcal{L}$, then
the density of probability of $\rho_i$ with respect to $\gamma_d$ is $\frac{d\rho_i}{\gamma_d}=f_i(2\pi)^{\frac{d}{2}}
e^{\frac{|x|_2^2}{2}}$. Write
\begin{eqnarray*}
	Ent_{\gamma_d}(\rho_i) &=&\int_{\mathbb{R}^d}f_i(x)(2\pi)^{\frac{d}{2}}e^{\frac{|x|_2^2}{2}}\log\left(f_i(x)
	(2\pi)^{\frac{d}{2}}e^{\frac{|x|_2^2}{2}}\right)d\gamma_d(x)\\
	&=&\int	f_i(x)(2\pi)^{\frac{d}{2}}e^{\frac{|x|_2^2}{2}}\log(f_i(x))d\gamma_d(x)
		+\int	f_i(x)(2\pi)^{\frac{d}{2}}e^{\frac{|x|_2^2}{2}}\log((2\pi)^{\frac{d}{2}})d\gamma_d(x)\\
	&~&
		+\int f_i(x)(2\pi)^{\frac{d}{2}}e^{\frac{|x|_2^2}{2}}\frac{|x|_2^2}{2}d\gamma_d(x)\\
	&=& Ent_{\mathcal{L}}(\rho_i)+\int \frac{1}{2}|x|_2^2d\rho_i(x)+\frac{d}{2}\log(2\pi).
\end{eqnarray*}
$\ast$ First term of the latter equality is relative entropy with respect to $\mathcal{L}$. Since $\|.\|^2$ is
strictly convex and differentiable, it suffices to follow \cite{AMB} to see that
it is convex along geodesics on $\mathcal{P}_2(\mathbb{R}^d,\|.\|)$.\\
$\ast$ Let us show that the second term of the latter equality $\mu \longmapsto \int \frac{1}{2}|x|_2^2d\mu(x)$
is $1-$convex along geodesics on $\mathcal{P}_2(\mathbb{R}^d,\|.\|)$. We know that the map $x\longmapsto \frac{1}{2}
|x|_2^2$ is $1-$convex i.e for all $x_1,~x_2 \in \mathbb{R}^d$ and $t$ in $[0,1]$:
\begin{eqnarray*}
	\frac{1}{2}|(1-t)x_1+tx_2|_2^2 &\leq& \frac{1-t}{2}|x_1|_2^2+\frac{t}{2}|x_2|_2^2-\frac{t(1-t)}{2}|x_1-x_2|_2^2\\
		&\leq& \frac{1-t}{2}|x_1|_2^2+\frac{t}{2}|x_2|_2^2-\frac{Kt(1-t)}{2}\|x_1-x_2\|^2.
\end{eqnarray*}
Consider an optimal coupling $\Pi$ (for $\|.\|^2$) between $\rho_0$ and $\rho_1$. Then integrating the previous inequality
w.r.t. $\Pi$, it comes:
\begin{eqnarray*}
	\int \frac{1}{2}|(1-t)x_1+tx_2|_2^2d\Pi(x_1,x_2)&\leq& \frac{1-t}{2}\int |x_1|_2^2d\rho_0(x_1)+\frac{t}{2}\int |x_2|_2^2
	d\rho_1(x_2) -\frac{Kt(1-t)}{2}W_{2}^2(\rho_0,\rho_1).
\end{eqnarray*}
$\ast$ Finally relative entropy with respect to $\gamma_d$ is $1-$convex along (constant speed) geodesics as sum of
$0-$convex and $1-$convex maps.
\end{demo}

We will apply Proposition \ref{EntRnNorm} in the following cases:
\begin{itemize}
	\item For the norm $|.|_p$ with $2\leq p < \infty$ ($K=1$).
	\item For the norm induced by $\|.\|_{k,\gamma}$ ($K=1/C_{k,\gamma}^2$).
\end{itemize}

We shall extend the result for the uniform norm $|.|_\infty$. This latter fact
is a priori not obvious since $|x-y|_\infty^2$ is neither strictly convex nor differentiable on $\mathbb{R}^d\backslash \{0\}$.\\

Now the question is: what is happening when $p$ equals $+\infty$ ? 
When one changes the cost function, two items change in the above inequality:
Wasserstein distances but also (constant speed) geodesics which
depend implicitely on an optimal transference plan, which depends itself on cost function.\\

Let $\rho_0$ and $\rho_1$ be definitively fixed, as probability measures on $\mathbb{R}^d$ with finite second moments (for $|.|_\infty$). In particular
$W_{2,p}(\rho_0,\rho_1)<\infty$ for all $p \geq 2$.
We know that for $p\geq 2$, there exists a unique optimal transference plan $\Pi_0^{(p)}$ between $\mu_0$ and $\mu_1$
(for cost function $c_p^2$). 
Thus we can watch behavior of the sequence $(\Pi_0^{(p)})_p$. When $p$ varies, cost function varies too.
In fact $(c_p)_p$ converges to $c$ and it would be interesting if the sequence of optimal coupling converge to an
optimal coupling for infinite cost.
Indeed it appears that up to a subsequence, $(\Pi_0^{(p)})_p$ weakly converges to
a probability measure which will be an optimal transference plan for the infinite cost. This fact combined with the property
of lower semicontinuity of the relative entropy, which one adds the nonincreasing of the following sequel
$$p\in \mathbb{N} \longmapsto W_{2,p}^2(\mu_0,\mu_1)$$
will yield $1-$convexity of relative entropy along geodesics on $(\mathbb{R}^d,|.|_\infty)$.\\

To prove the weak convergence of $(\Pi_0^{(p)})_p$ a first (easy but useful) remark is that the sequel is tight.
This yields, thanks to Prokohov's Theorem, there exists a subsequence $(\Pi_0^{(p_k)})_{p_k}$ that we will denote 
by $(\Pi_0^{(p)})_p$ again, converging weakly to a measure $\Pi^{\infty}$. It is easy to check that $\Pi^{\infty}$ is a coupling of
$\mu_0$ and $\mu_1$. And we wish to see that $\Pi^{\infty}$ is optimal for the cost $c^2$, where $c(x,y):=|x-y|_{\infty}$.
The difficulty lies in the fact that $(c_p)_p$ converges to $c$ but \emph{not uniformly}. Nevertheless, it will work because 
the convergence is uniform on all compact subsets of $\mathbb{R}^d$.\\

The next Lemma which appears to be an essential point for this purpose, is taken from \cite{AMB}. Recall that for a probability measure $\mu$, its \emph{support} $Supp(\mu)$, is defined as the smallest closed set on which $\mu$ is concentrated (i.e. $\mathbb{R}^d\backslash Supp(\mu)$ is $\mu-$negligible).

\begin{lemme}\label{supportlimite}
For all $x \in Supp(\Pi^{\infty})$, there exist $x_p \in Supp(\Pi_0^{(p)})$ such that $lim (x_p)=x$.
\end{lemme}
	
For self contained paper, here we include the proof.

\begin{demo}
Let $x \in Supp(\Pi^{\infty}) \subset \mathbb{R}^d\times \mathbb{R}^d$. Thus $\Pi^{\infty}$ attributes a nonzero mass for all open ball of center $x$. So let $k \in \mathbb{N}^{\star}$, we have $\Pi^{\infty}(B(x,1/k))>0$, where $B(x,1/k)$ is the open ball centered at $x$
and radius $1/k$. Thanks to the weakly convergence, we have:
$$\liminf_{p\longrightarrow +\infty} \Pi_0^{(p)}(B(x,1/k)) \geq \Pi^{\infty}(B(x,1/k))>0.$$
This inequality let us define an increasing sequence $(j_k)_k$ such that: $j_0:=0$ and for $k>0$
$$j_k := \min \{p\in \mathbb{N},~p>j_{k-1},~\forall n\geq p~: Supp(\Pi_0^{(n)})\cap B(x,1/k) \neq \varnothing\}.$$
The increase yields that for all $k \in \mathbb{N}$, there exists $j_k \leq p < j_{k+1}$ such that we can pick up a point
$x_p \in Supp(\Pi_0^{(p)})\cap B(x,1/k)$. The sequence $(x_p)_p$ converges to $x$ and for all $p\in \mathbb{N}$, we have $x_p \in Supp(\Pi_0^{(p)})$.
\end{demo}

\begin{rmq}
In fact this proof (hence this lemma) is true for all sequence of measures converging weakly to another measure, on any metric space.
\end{rmq}

\begin{prop}
$\Pi^{\infty}$ is optimal for the cost $c^2$.
\end{prop}

Here we use the equivalence between \emph{optimality for $c^2$} and \emph{$c^2$-cyclical monotonicity} (see \cite{VIL} for
continuous and real valued costs $c$).

\begin{demo}
It is then sufficient to prove that $Supp(\Pi^{\infty})$ is $c-$cyclically monotone.
Let $N \in \mathbb{N}^{\star}$ and $(x_1,y_1),\dots,(x_N,y_N) \in Supp(\Pi^{\infty})$. Since $(\Pi_0^{(p)})_p$ converges weakly
to $\Pi^{\infty}$, we can applicate the Lemma \ref{supportlimite}: for all $i=1,\dots N$, there exists 
$(x_i^p,y_i^p) \in Supp(\Pi_0^{(p)})$ such that $\lim (x_i^p,y_i^p)=(x_i,y_i)$. Thus $(x_1^p,y_1^p),\dots,(x_N^p,y_N^p)\in Supp
(\Pi_0^{(p)})$ which is $c_p^2-$cyclically monotone, since $\Pi_0^{(p)}$ is optimal for the cost $c_p$. Then the inequality 
\begin{equation}\label{cpcycle}
	\sum_{i=1}^Nc_p^2(x_i^p,y_i^p)\leq \sum_{i=1}^Nc_p^2(x_i^p,y_{i+1}^p)
\end{equation}
holds, with $y_{N+1}:=y_1$. And it is elementary to check that the sets 
\begin{eqnarray*}
\cup_{p\geq 2}{ \{ (x_1^p,y_1^p),\dots,(x_N^p,y_N^p)\} \bigcup \{(x_1,y_1),\dots,(x_N,y_N)\}},\\
\cup_{p\geq 2}{ \{ (x_1^p,y_2^p),\dots,(x_N^p,y_1^p)\}\bigcup \{(x_1,y_2),\dots,(x_N,y_1)\}},
\end{eqnarray*}
are compact of $\mathbb{R}^d\times \mathbb{R}^d$.
But since $(c_p)_p$ converges uniformly on compact subsets of $\mathbb{R}^d\times \mathbb{R}^d$ to $c$, we get from (\ref{cpcycle}), taking
the limit with $p\rightarrow +\infty$:
$$\sum_{i=1}^Nc^2(x_i,y_i)\leq \sum_{i=1}^Nc^2(x_i,y_{i+1}).$$
That is exactly the definition of $c^2-$cyclically monotone for $Supp(\Pi^{\infty})$.
\end{demo}

Because of non strict convexity of $|.|_\infty$, $(\mathbb{R}^d,|.|_\infty)$ is
a branching space: there exists many constant speed geodesics between two probability measures.
Finally one can conclude this section with the following result:

\begin{prop}\label{EntRnInfini}
	Relative convexity w.r.t. $\gamma_d$ on $(\mathbb{R}^d,|.|_\infty)$ is \underline{weakly} $1-$convex along 
	(constant speed) geodesics on $(\mathcal{P}_2(\mathbb{R}^d),W_{2,\infty})$.
\end{prop}

\begin{demo}
For $p \in [2,+\infty)$ we consider optimal transfere plans $\Pi_0^{(p)} \in \Gamma_0^{p}(\rho_0,\rho_1)$ and constant speed geodesics
$\rho_t^{(p)}$ associated. Applying Proposition \ref{EntRnNorm} with $|.|_p$ norms, we get:
\begin{equation}\label{EntRn0}
	Ent_{\gamma_d}(\rho_t^{(p)})\leq (1-t)Ent_{\gamma_d}(\rho_0)+tEnt_{\gamma_d}(\rho_1)-\frac{t(1-t)}{2}W_{2,p}^2(\rho_0,\rho_1),
\end{equation}
for all $p\geq 2$. But for $x\in \mathbb{R}^d$ and for $p\geq 2$, $|x|_p \geq |x|_\infty$. Then:
\begin{eqnarray*}
	\int_{\mathbb{R}^d\times \mathbb{R}^d} |x-y|_p^2d\Pi_0^{(p)}(x,y)&\geq&\int_{\mathbb{R}^d\times \mathbb{R}^d} |x-y|_\infty ^2d\Pi_0^{(p)}(x,y),\\
	W_{2,p}^2(\rho_0,\rho_1) &\geq & \inf_{\Pi\in \Gamma(\mu_0,\mu_1)}\int_{\mathbb{R}^d\times \mathbb{R}^d} |x-y|_\infty ^2d\Pi(x,y) = W_{2,\infty}^2(\rho_0,\rho_1),\\
	\Longrightarrow \liminf_p W_{2,p}^2(\rho_0,\rho_1) &\geq & W_{2,\infty}^2(\rho_0,\rho_1).
\end{eqnarray*}
Moreover, the sequel $(\Pi_0^{(p)})_p$ weakly converges to $\Pi^{\infty}$ and this coupling is optimal
for $|.|_\infty^2$ thanks to the previous proposition. Hence the sequel $(\rho_t^{(p)})_p$ weakly converges
to $\rho_t^{\infty}$ for all $t \in [0,1]$. But since the relative entropy is lower semi-continuous, we have:
$$\liminf_p Ent_{\gamma_d}(\rho_t^{(p)}) \geq Ent_{\gamma_d}(\rho_t^{\infty}).$$
Finally, combining this two arguments, taking the liminf in the inequality (\ref{EntRn0}) with respect to $p$,
$$Ent_{\gamma_d}(\rho_t^{\infty})\leq (1-t)Ent_{\gamma_d}(\rho_0)+tEnt_{\gamma_d}(\rho_1)-\frac{t(1-t)}{2}W_{2,\infty}^2(\rho_0,\rho_1).$$
\end{demo}

\subsection{Wiener space case}

Consider the same notations of the previous settings and denote by $\mathcal{A}_n$ the sub $\sigma-$field
on $W$ generated by $\pi_n$. Let us recall the following results:
\begin{enumerate}
	\item $(\pi_n)_\#\mu = \gamma_n$ the standart Gaussian measure on $V_n$.
	\item For $i=0,1$, $\rho_i^n:=(\pi_n)_\# \rho_i$ is absolutely continuous with respect to $\gamma_n$ with density $f_i^n$ where
	$f_i^n\circ \pi_n =\mathbb{E}[f_i|\mathcal{A}_n] =:\hat{f}_i^n$.
	\item For $i=0,1$, $\hat{f}_i^n\longrightarrow f_i$ in $L^1(W,\mu)$.
\end{enumerate}

\begin{rmq}
	In the following proof, we need to apply Proposition \ref{EntRnInfini} which deals with the infinite norm
	in $\mathbb{R}^n$. It turns out that the above approximation by $(\pi_n)_n$ doesn't lead us to the infinite norm.
	We get only a modified (by an inversible matrice) infinite norm. But the previous section can be applied to
	this modified norm, so that Proposition \ref{EntRnInfini} also holds in our case.
\end{rmq}

\begin{thm}\label{EntW}
	Relative entropy w.r.t. $\mu$ on $(W,|.|_\infty)$ is \underline{weakly} $1-$convex along (constant speed) geodesic
	on $(\mathcal{P}_2(W),W_{2,\infty})$.
\end{thm}

\begin{demo}
\textbf{Step 1.}
Since $\hat{f}_i^n=f_i^n\circ \pi_n \longrightarrow f_i$, we have: 
$$\hat{\rho}_i^n:= (f_i^n\circ \pi_n)\mu \stackrel{weakly}{\longrightarrow} f_i\mu = \rho_i~~~~ i=0,1.$$
Then in \cite{AMB}, it is proved that for $\hat{\Pi}_n^{\infty} \in \Gamma_0^{\infty}(\hat{\rho}_0^n,\hat{\rho}_1^n)$ there exists
a subsequence of $(\hat{\Pi}_n^{\infty})_n$ denoted by $(\hat{\Pi}_{n_k}^{\infty})_k$ which weakly converges to $\hat{\Pi}^{\infty}$ and we
know this limit point belongs to $\Gamma_0^{\infty}(\rho_0,\rho_1)$.\\

\vspace{5pt}
\noindent \textbf{Step 2.} $W_{2,\infty}(\hat{\rho_0}^n,\hat{\rho_1}^n)\leq W_{2,\infty}(\rho_0^n,\rho_1^n)$\\
We have the decomposition $W=V_n \oplus V_n^{\bot}$. Then define a probability measure $\Pi \in \Gamma(\hat{\rho}_0^n,
\hat{\rho}_1^n)$ as
$$\int_{W\times W}\varphi (x,y)d\Pi(x,y) := \int_{V_n^{\bot}}\int_{V_n\times V_n}\varphi(x_1+z,x_2+z)d\Pi_n(x_1,x_2)d\rho_0^{\bot}(z),$$
for $\varphi$ bounded continuous function, and $\Pi_n \in \Gamma_0^{\infty}(\rho_0^n,\rho_1^n)$.
So:
\begin{eqnarray*}
	\int_{W\times W}|x-y|_\infty^2d\Pi(x,y)&=&\int_{V_n\times V_n}|x_1-x_2|_\infty^2d\Pi_n(x_1,y_1)\\
		&=&W_{2,\infty}^2(\rho_0^n,\rho_1^n),\\
		\Longrightarrow W_{2,\infty}^2(\hat{\rho}_0^n,\hat{\rho}_1^n) &\leq & W_{2,\infty}^2(\rho_0^n,\rho_1^n).
\end{eqnarray*}

\noindent \textbf{Step 3.} \\
Now we have $\Pi_{n_k}^\infty :=(\pi_{n_k},\pi_{n_k})_\# \hat{\Pi}_{n_k}^\infty \in \Gamma_0^\infty(\rho_0^{n_k},\rho_1^{n_k})$ and
if we consider the associated constant speed geodesic $\rho_t^{n_k}$, we can apply Proposition \ref{EntRnInfini} to obtain:
\begin{equation}\label{EntRn2}
	Ent_{\gamma_{n_k}}(\rho_t^{n_k})\leq (1-t)Ent_{\gamma_{n_k}}(\rho_0^{n_k})+tEnt_{\gamma_{n_k}}(\rho_1^{n_k})-\frac{t(1-t)}{2}
	W_{2,\infty}^2(\rho_0^{n_k},\rho_1^{n_k})~~~~\forall t\in [0,1].
\end{equation}
Let $t\in [0,1] \longmapsto \hat{\rho}_t^{n_k}:=((1-t)Pr_1+tPr_2)_\#\hat{\Pi}_{n_k}^\infty$ be the constant speed geodesic associated
to $\hat{\Pi}_{n_k}^\infty$. We have:
\begin{eqnarray*}
	Ent_{\gamma_{n_k}}(\rho_0^{n_k}) &=& \int_{V_{n_k}}f_0^{n_k}(x)\log(f_0^{n_k}(x))d\gamma_{n_k}(x)\\
		&=& \int_W f_0^{n_k}(\pi_{n_k}(y))\log(f_0^{n_k}(\pi_{n_k}(y)))d\mu(y)\\
		&=& \int_W \hat{f}_0^{n_k}(y)\log(\hat{f}_0^{n_k}(y))d\mu(y) = Ent_\mu(\hat{\rho}_0^{n_k}).
\end{eqnarray*}
And on the same way, we can easily prove this for $\rho_1^{n_k}$ and $\rho_t^{n_k}$ where $t \in (0,1)$. Thus (\ref{EntRn2}) and 
Step 2. yield
\begin{equation}\label{EntRn3}
	Ent_{\mu}(\hat{\rho}_t^{n_k})\leq (1-t)Ent_{\mu}(\hat{\rho}_0^{n_k})+tEnt_{\mu}(\hat{\rho}_1^{n_k})-\frac{t(1-t)}{2}W_{2,\infty}^2
	(\hat{\rho}_0^{n_k},\hat{\rho}_1^{n_k})~~~~\forall t\in [0,1].
\end{equation}

\noindent \textbf{Step 4.} \\
Thanks to Step 1., we have:
\begin{eqnarray*}
	W_{2,\infty}^2(\rho_0,\rho_1) &=& \int_{W\times W}|x-y|_\infty^2d\hat{\Pi}^\infty(x,y) \\
		&\leq & \liminf_{k} \int_{W\times W} |x-y|_\infty^2d\hat{\Pi}_{n_k}^\infty (x,y) = \liminf_k W_{2,\infty}^2(\hat{\rho}_0^{n_k},
		\hat{\rho}_1^{n_k}).
\end{eqnarray*}
Hence for $\epsilon >0$, there exists $N \in \mathbb{N}$ such that:
$$W_{2,\infty}^2(\hat{\rho}_0^{n_k},\hat{\rho}_1^{n_k})+\epsilon \geq W_{2,\infty}^2(\rho_0,\rho_1)~~~~\forall k \geq N.$$
Now fix $k \geq N$. By the Jensen's inequality, we have for $i=0,1$:
\begin{eqnarray*}
	Ent_\mu(\hat{\rho}_i^{n_k}) &=&\int_{W}\hat{f}_i^{n_k}(y)\log(\hat{f}_i^{n_k}(y))d\mu(y) \\
		&=& \int_W \mathbb{E}[f_i|\mathcal{A}_{n_k}](y)\log(\mathbb{E}[f_i|\mathcal{A}_{n_k}](y))d\mu(y)\\
		&\leq & \int_W \mathbb{E}[f_i\log(f_i)|\mathcal{A}_{n_k}](y)d\mu(y)\\
		&\leq & \int_W f_i(y)\log(f_i(y))d\mu(y) = Ent_\mu(\rho_i),
\end{eqnarray*}
and then:
\begin{eqnarray*}
	Ent_{\mu}(\hat{\rho}_t^{n_k}) &\leq & (1-t)Ent_\mu(\hat{\rho}_0^{n_k})+tEnt_\mu(\hat{\rho}_1^{n_k})-\left(\frac{t(1-t)}{2}W_{2,\infty}
	^2(\rho_0,\rho_1)-\epsilon\right)\\
	&\leq & (1-t)Ent_\mu(\rho_0)+tEnt_\mu(\rho_1)-\left(\frac{t(1-t)}{2}W_{2,\infty}^2(\rho_0,\rho_1)-\epsilon\right).
\end{eqnarray*}
Moreover:
\begin{eqnarray*}
	\mathrm{Since}~~~~\hat{\Pi}_{n_k}^\infty &\stackrel{weakly}{\longrightarrow}& \hat{\Pi}^\infty~
	\mathrm{which~is~optimal~due~to~Lemma~3.1,}\\
	\mathrm{we~have}~~~~\hat{\rho}_t^{n_k} &\stackrel{weakly}{\longrightarrow}& \rho_t:=((1-t)Pr_1+tPr_2)_\# \hat{\Pi}^\infty.
\end{eqnarray*}
Denote by $R$ the right hand side of the previous inegality, and by compacity of $\{\nu =\rho \mu
\in \mathcal{P}_2(W),~Ent_\mu(\nu)\leq R\}$
with respect to the weak topology, we have:
$$Ent_\mu(\rho_t)\leq (1-t)Ent_\mu(\rho_0)+tEnt_\mu(\rho_1)-\left(\frac{t(1-t)}{2}W_{2,\infty}^2(\rho_0,\rho_1)-\epsilon\right).$$
Finally letting $\epsilon \longrightarrow 0$ we get:
$$Ent_\mu(\rho_t)\leq (1-t)Ent_\mu(\rho_0)+tEnt_\mu(\rho_1)-\frac{t(1-t)}{2}W_{2,\infty}^2(\rho_0,\rho_1)~~~~\forall t\in[0,1].$$

\end{demo}

The same proof holds for the cost $\|.\|_{k,\gamma}$ by applying in Step 3. 
Proposition \ref{EntRnNorm} with $\|.\|_{k,\gamma}\leq C_{k,\gamma}|.|_2$
instead of Proposition \ref{EntRnInfini}.
Even better, since the latter norm is strictly convex, $(W,\|.\|_{k,\gamma})$
is a non branching space. Hence we have the Theorem:

\begin{thm}
	Relative entropy w.r.t. $\mu$ on $(W,\|.\|_{k,\gamma})$ is $1/C_{k,\gamma}^2-$convex along (constant speed) geodesics
	on $(\mathcal{P}_2(W),W_{2,(k,\gamma)})$.
\end{thm}

\begin{itemize}
	\item In Step 2. of first proof we have in fact equality $W_{2,\infty}(\hat{\rho}_0^n,\hat{\rho}_1^n)=W_{2,\infty}(\rho_0^n,\rho_1^n)$. This is provided
by the fact that $|\pi_n(w)|_\infty\leq |w|_\infty$ for all $w\in W$. Indeed
we know that $\hat{\Pi}_n^{\infty}$ is optimal. Define $\Pi_n:=(\pi_n,\pi_n)_\# \hat{\Pi}_n^{\infty} \in \mathcal{P}(V_n\times V_n)$. Since
the norm on $V_n$ is less than the norm on $W$, we have:
\begin{eqnarray*}
	\int_{V_n\times V_n} |x-y|_\infty ^2d\Pi_n(x,y) &=& \int_{W\times W}|\pi_n(x)-\pi_n(y)|_\infty^2d\hat{\Pi}_n^{\infty}(x,y)\\
		&\leq & \int_{W\times W}|x-y|_\infty ^2d\hat{\Pi}_n^{\infty}(x,y)= W_{2,\infty}^2(\hat{\rho}_0^n,\hat{\rho_1}_1^n),\\
		\Longrightarrow W_{2,\infty}^2(\rho_0^n,\rho_1^n) &\leq & W_{2,\infty}^2(\hat{\rho}_0^n,\hat{\rho}_1^n).
\end{eqnarray*}
	\item In the second norm, it is still not clear if $\|\pi_n(w)\|_{k,\gamma}\leq \|w\|_{k,\gamma}$ for any $w\in W$.
\end{itemize}

\begin{appendices}

\begin{thm}\emph{\textbf{\texttt{Rademacher's Theorem.}}} If $\varphi \in L^p(W)$ is a $\|.\|_{k,\gamma}-$convex 
map then $\varphi$ belongs to $\mathbb{D}_1^p(W)$.\end{thm}

\begin{demo}
Our discussion above gives us for $\mu-$almost all $x$ and $y$ in $W$:
$$|\varphi(x)-\varphi(y)|\leq c(x,y)=\|x-y\|_{k,\gamma}\leq 2|x-y|_H.$$
This can be rewritten for all $w\in W$ and all $h \in H$ by:
$$|\varphi(w+h)-\varphi(w)|\leq 2 |h|_H.$$
Fix $l\in W^{\star}\subset H$. Consider the set
$$\Lambda(l)=\left\{w\in W,~G(w,l):=\lim_{\epsilon \rightarrow 0} \frac{\varphi(w+\epsilon l)-\varphi(w)}{\epsilon}~\mathrm{exists}\right\}$$
and an orthornomal basis $(y_n)_n$ of $Y$, where $Y$ is such that $W=Y\oplus Span(l)$. By denstiy one can write:
\begin{eqnarray*}
	\mu\left(\Lambda(l)^c\right) &=& \mu\left(\left\{y_n+tl \in W;~\lim_{\epsilon \rightarrow 0} \frac{\varphi(y_n+(\epsilon+t)l )-\varphi(y_n+tl)}{\epsilon}~\mathrm{not~exists}\right\}\right)\\
		&=&Leb\left(\bigcup_n \left\{t \in \mathbb{R};~G(y_n+tl,l)~\mathrm{not~exists}\right\}\right)\\
		&\leq&\sum_n Leb \left(\left\{t \in \mathbb{R};~G(y_n+tl,l)~\mathrm{not~exists}\right\}\right).
\end{eqnarray*}
And the measure of $\left\{t \in \mathbb{R};~G(y_n+tl,l)~\mathrm{does~not~exist}\right\}$ is equal to zero, thanks to the Rademacher's Theorem on $\mathbb{R}$.
Then $\mu(\Lambda(l))=1$ and from now we extend $G(.,l)$ on $W$ by set $G(w,l)=0$ for all $w\notin \Lambda(l)$.
By assumption we have for all $\epsilon >0$ and $w\in W$:
$$	\frac{|\varphi (w+\epsilon l)-\varphi (w)|}{\epsilon} \leq  2 \frac{\epsilon |l|_H}{\epsilon} =2|l|_H$$
This implies that:
$$|G(w,l)|\leq C|l|_H~~~~\forall w\in W,~\mathrm{and}~\forall l \in W^\star.$$
Now consider:
$$\mathcal{A}_n:=\left\{\sum_{k=1}^n a_ke_k;~a_k\in \mathbb{Q}~\mathrm{and}~\sum_{k=1}^na_k^2=1\right\},~~\mathcal{A}=\bigcup_{n\geq 1}\mathcal{A}_n.$$
If $l=\sum_{k=1}^na_ke_k \in \mathcal{A}$ then we have $\mu-$almost surely:
$$G(w,l)=\sum_{k=1}^n(l,e_k)_HG(w,e_k).$$
Therefore for $l\in \mathcal{A}$, it exists $n\in \mathbb{N}$ such that $l$ belongs to $\mathcal{A}_n$ and a Borel subset
$B_n$ of $W$ with plain $\mu-$measure, such that we have:
$$G(w,l)=\sum_{k=1}^n(l,e_k)_HG(w,e_k)~~~~\forall w\in B_n.$$
Since $\mathcal{A}$ is countable, take $B=\cap_n B_n$ which has plain $\mu-$measure and satisfies:
$$G(w,l)=\sum_{k=1}^\infty (l,e_k)_HG(w,e_k)~~~~\forall w \in B, ~\forall l \in \mathcal{A}.$$
Now define the map $G(w)=\sum_{k=1}^\infty G(w,e_k)e_k$ for $w\in B$ and equals to $0$ if $w\notin B$.
Notice that:
$$(G(w),l)_H=\sum_{k\geq 1}G(w,e_k)(e_k,l)_H=G(w,l)~~~~\forall w\in B,~\forall l\in \mathcal{A}.$$
Then we have for $w\in B$:
\begin{eqnarray*}
	|G(w)|_H^2 & = & \left(G(w),\sum_{k\geq 1}G(w,e_k)e_k \right)_H\\
		&=& G\left(w,\sum_{k\geq 1}G(w,e_k)e_k\right),
\end{eqnarray*}
and the latter term being less than $2|G(w)|_H$, we get a bound from above of $|G(w)|_H$:
\begin{equation}\label{Gmaj}
	|G(w)|_H\leq 2~~~~\forall w\in B.
\end{equation}

At last we have for $h\in H$:
\begin{eqnarray*}
	\left\|\frac{\varphi(w+\epsilon h)-\varphi(w)}{\epsilon}-(G(w),h)_H\right\|_p^p &\leq &
		\int_{\Lambda(h)} \left| \frac{\varphi(w+\epsilon h)-\varphi(w)}{\epsilon}-(G(w),h)_H\right|^pd\mu(w).
\end{eqnarray*}
But thanks to (\ref{Gmaj}) and Cauchy-Schwartz inequality, we have:
$$\left| \frac{\varphi(w+\epsilon h)-\varphi(w)}{\epsilon}-(G(w),h)_H\right|\leq 2^2|h|_H$$
and this latter function belongs to $L^p(W)$. So applying dominated convergence theorem, we get
$$\lim_{\epsilon \rightarrow 0}	\left\|\frac{\varphi(w+\epsilon h)-\varphi(w)}{\epsilon}-(G(w),h)_H\right\|_p^p=0.$$
This shows that $D_h\varphi(w)=(G(w),h)_H$ exists $\mu-$almost surely in $w$. Moreover 
$G(w)$ is the gradient of $\varphi$ at $w$. Finally
$G$ belongs to $L^p(W,H)$ i.e. $\varphi \in \mathbb{D}_1^p(W)$.
\end{demo}

\begin{defin}
	We say that a function $f:W\longrightarrow W$ is \emph{locally $\|.\|_{k,\gamma}-$Lipschitz} if for all $R>0$,
	there exists $C_R>0$ such that:
	$$|f(w)-f(w')|\leq C_R \|w-w'\|_{k,\gamma}~~~~\forall w,w'\in B_{k,\gamma}(R):=\{w\in W;~\|w\|_{k,\gamma}\leq R\}.$$
\end{defin}

\begin{prop}\label{SupConvex} Let $(f_i)_{i\in I}$ be a uniform family of locally $\|.\|_{k,\gamma}-$Lipschitz real-valued 
functions defined on an open subset $U$ of $W$. If the function
$$f(x):=\sup_{i\in I}f_i(x)$$
is finite everywhere then $f:U\longrightarrow \mathbb{R}$ is also locally $\|.\|_{k,\gamma}-$Lipschitz.
\end{prop}

\begin{demo}
	Let $R>0$ such that $B_{k,\gamma}(R)\subset U$. There exists $C_R>0$ (not depending on $i$) such that
	for all $w, w'\in B_{k,\gamma}(R)$:
\begin{eqnarray*}
	|f_i(w)-f_i(w')|&\leq& C_R\|w-w'\|_{k,\gamma}~~~~\forall i\in I.
\end{eqnarray*}
	For each $w \in B_{k,\gamma}(R)$ by definition of $f$ there is a sequence 
	$(i_n)_n$ (depending on $w$) such that $\lim_{n}f_{i_n}(w)=f(w)$. Moreover
	$f_{i_n}(w')\leq f(w')$ for all $w'\in B_{k,\gamma}(R)$ then
$$f_{i_n}(w)-f(w')\leq f_{i_n}(w)-f_{i_n}(w')\leq C_R \|w-w'\|_{k,\gamma}.$$
	Passing to the limit in the previous inequality, we get:
	$$f(w)-f(w')\leq C_R\|w-w'\|_{k,\gamma}.$$
	Exchanging $w$ and $w'$, this yields $f$ is locally $\|.\|_{k,\gamma}-$Lipschitz.
\end{demo}

\end{appendices}

\thanks{\textbf{Acknowledgements.} I would like to thank Professor Shizan Fang, for his
suggestion, and so Nicolas Juillet for useful discussions about
general facts of optimal transportation.}

\nocite{*}
\bibliographystyle{plain}
\bibliography{biblio}

\begin{thebibliography}{10}

\bibitem{AIRMALL}
H.~Airault and P.~Malliavin.
\newblock {Integration geometrique sur l'espace de Wiener}.
\newblock {\em Bulletin des Sciences Mathematiques}, pages 3--52, 1988.

\bibitem{AMB0}
L.~Ambrosio.
\newblock {Optimal transport maps in Monge-Kantorovich problem}.
\newblock {\em Proceedings of the International Congress of Mathematicians,
  Vol. III}, pages 131--140, 2002.

\bibitem{USERG}
L.~Ambrosio and N.~Gigli.
\newblock A user's guide to optimal transport.
\newblock 2011.

\bibitem{AMB}
L.~Ambrosio, N.~Gigli, and G.~Savare.
\newblock {\em {Gradient Flows in Metric Spaces and in the Space of Probability
  Measures}}.
\newblock Lectures in Mathematics.

\bibitem{AMB2}
L.~Ambrosio, B.~Kirchheim, and A.~Pratelli.
\newblock Existence of optimal transport maps for crystalline norms.
\newblock {\em Duke Mathematical Journal}, pages 207--241, 2004.

\bibitem{BIANCAV}
S.~Bianchini and F.~Cavalletti.
\newblock The monge problem for distance cost in geodesic spaces.
\newblock {\em Submitted Paper}, 2009.

\bibitem{CARAV}
L.~Caravenna.
\newblock {An existence result of the Monge problem in $\mathbb{R}^n$ with norm
  cost functions}.
\newblock 2010.

\bibitem{CAVA}
F.~Cavalletti.
\newblock {The Monge Problem in Wiener space}.
\newblock {\em Calculus of Variations}, 2011.

\bibitem{CHAMP}
T.~Champion and L.~De~Pascale.
\newblock {The Monge problem in $\mathbb{R}^d$}.
\newblock {\em Duke Mathematical Journal}, 2010.

\bibitem{ENCHSTRO}
O.~Enchev and W.~Stroock.
\newblock Rademacher's theorem for wiener functionals.
\newblock {\em The Annals of Probability}, pages 25--33, 1993.

\bibitem{FANG0}
S.~Fang.
\newblock {\em {Introduction to Malliavin Calculus}}.
\newblock Mathematics Series for Graduate Students, 2003.

\bibitem{FANGNOL}
S.~Fang and V.~Nolot.
\newblock {Sobolev estimates for optimal transport maps on Gaussian spaces}.
\newblock {\em arXiv}, 2012.

\bibitem{FANG}
S.~Fang, J.~Shao, and K-T. Sturm.
\newblock {Wasserstein space over the Wiener space}.
\newblock {\em Probab. Theory Related Fields}, pages 535--565, 2010.

\bibitem{FAFI}
A.~Fathi and F.~Figalli.
\newblock Optimal transportation on non-compact manifolds.
\newblock {\em Israel J. Math.}, pages 1--59, 2010.

\bibitem{FEY}
D.~Feyel and A.S. Ustunel.
\newblock {Monge-Kantorovitch measure transportation and Monge-Ampere equation
  on Wiener space}.
\newblock {\em Probab. Theory Related Fields}, pages 347--385, 2004.

\bibitem{GANGMC}
W.~Gangbo and R.J. McCann.
\newblock The geometry of optimal transportation.
\newblock {\em ActaMath.}, pages 113--161, 1996.

\bibitem{GANGOL}
W.~Gangbo and V.~Oliker.
\newblock Existence of optimal maps in the reflector-type problems.
\newblock {\em ESAIM Control Optim. Calc. Var.}, pages 93--106, 2007.

\bibitem{Lott}
J.~Lott and C.~Villani.
\newblock Ricci curvature for metric-measure spaces via optimal transport.
\newblock {\em Annals of Mathematics}, pages 903--991, 2009.

\bibitem{McC}
R.J. McCann.
\newblock Existence and uniqueness of monotone measure-preserving maps.
\newblock {\em Duke Mathematical Journal}, pages 309--323, 1995.

\bibitem{STURM}
K.T. Sturm.
\newblock {On the Geometry of Metric Measure Spaces I.}
\newblock {\em Acta Math}, pages 65--131, 2006.

\bibitem{VIL}
C.~Villani.
\newblock {\em Optimal transport, old and new}.
\newblock Grundlehren der mathematischen Wissenschaften, 2009.

\end{thebibliography}

\end{document}